\documentclass{article}
 \UseRawInputEncoding
\usepackage{amsmath}
\font\tengoth=eufm10
\font\sevengoth=eufm7
\font\fivegoth=eufm5
\newfam\gothfam
\textfont\gothfam=\tengoth \scriptfont\gothfam=\sevengoth
\scriptscriptfont\gothfam=\fivegoth

\def\blacksquare{\hbox to .60em{\vrule width .60em height .60em}}

  \font\bb=msbm10 
\catcode`\@=12  

\def\é{\'e}
\def\{\`e}
\def\?{\`a}
\def\{\`u}
\def\{\c c}
\def\hb {\hfil \break}
\def\n {\vskip 0.2cm \noindent }
\def\scirc{\,{\raise 0.8pt\hbox{$\scriptstyle\circ$}}\,}
\def\ins{\,{\raise 0.2cm \hbox{ $\scriptstyle \circ$}}\,}

\def  \é{\'e}
\def\è{\`e}
\def\à{\`a}
\def\ù{\`u}
\def\ç{\c c$\!\!\!$}

 \date{ }

\begin{document}
  \centerline {\bf Relations ab\éliennes des tissus ordinaires de codimension arbitraire }
\medskip


 
 
 
\title{Relations ab\éliennes des tissus ordinaires de codimension arbitraire}
\author{  Daniel Lehmann}

 \bigskip
 \centerline{  Daniel Lehmann}
 
 
\n {\bf Abstract}
 
 {\it  We  generalize to  webs of any codimension results   already known in codimension one. 
 
Given a holomorphic $d$-web $\cal W$ of codimension $q$ $(q\leq n-1)$ in an ambiant $n$-dimensional holomorphic manifold $U$, we define for any integer  $p$  $(1\leq p\leq q)$ the condition for such a web to be \emph{$p$-ordinary} $($resp. \emph{strongly $p$-ordinary}$)$. If this condition is satisfied, we then prove  that its  $p$-rank $r_p({\cal W})$ $\bigl($resp. its closed $p$-rank $\widetilde r_p({\cal W})\bigr)$,  i.e. the maximal dimension of the vector space of the germs of $p$-abelian relations $($resp. of closed $p$-abelian relations$)$ at a point $m$ of $U$, 
 is finite.  We then give  an upper-bound $\pi_p^0(n,d,q)$ $\bigl($resp. $\pi'_p(n,d,q)\bigr)$ for  these ranks. 
 
 Moreover, for some values of $d$, and we then say then that the web is  \emph{$p$-calibrated} $($resp. \emph{strongly $p$-calibrated}$)$, we define a tautological holomorphic connection on a holomorphic vector bundle of rank $\pi_p^0(n,d,q)$ $\bigl($resp. $\pi'_p(n,d,q)\bigr)$, for which the sections with  vanishing  covariant derivative  may be identified with $p$-abelian relations $($resp.  closed $p$-abelian relations$)$. The curvature of this connection is then an obstruction for the rank $r_p({\cal W})$ $\bigl($resp. $\widetilde r_p({\cal W})\bigr)$  to be maximal. 
 
The main change in this new version  is the correction of  a mistake  $($proposition 4, section 
   6-5$)$ of the first one  :   the 1-rank of the concerned web is not 0 as we claimed, but 1.  However, the important corollary remains true : even at the level of germs, some  2-abelian relation  exhibited by Goldberg in $ [G]$ on some web of codimension 2 in an ambiant space of dimension 4,
  is   the coboundary of none 1-abelian relation.   The section 7, devoted to this correction,  is self content, not depending  on  the previous results of the paper.}

\vskip 2cm
 
 \n  { \bf Table des mati\ères :}
 \bigskip
 
  1- Introduction. \bigskip
  
  2- D\éfinition des $p$-relations ab\éliennes. Notations et calculs pr\éliminaires. \bigskip
  
  3- Relations ab\éliennes des tissus ordinaires.\bigskip
  
  4-  Relations ab\éliennes ferm\ées.\bigskip
  
  5-  Tissus $p$-calibr\és et connexions.\bigskip
  
  6- Exemples. \bigskip
  
  7- Un germe de   relation ab\élienne ferm\ée qui n'est pas un cobord. \bigskip
  
  R\éf\érences. 
  
    \vskip 1.5cm
 \n {\bf Mots cl\és :} tissus (fortement) $p$-ordinaires,    $p$-relations ab\éliennes (ferm\ées), $p$-rang (ferm\é), connexion, courbure. 
 
\n {\bf  Classification AMS :} 53A60 (14C21, 53C05, 14H45)

 \pagebreak
  \section{Introduction} 
  
Soit $\cal W$ un $d$-tissu holomorphe  de codimension $q$, compl\ètement d\écomposable,  dans une vari\ét\é holomorphe $U$ de dimension $n$ ($1\leq q\leq n-1$), d\éfini par $d$ feuilletages $({\cal F}_i)_{1\leq i\leq d}$ de codimension $q$ en position g\én\érale.  Rappelons ([Gr]) qu'une $p$-relation ab\élienne ($1\leq p\leq q$) est la donn\ée d'une famille $(\omega_i)_{1\leq i\leq d}$ de $p$-formes diff\érentielles sur $U$, o\ù  chaque forme $\omega_i$ est    ${\cal F}_i$-basique,   et $\sum_i \omega_i=0$,. L'ensemble $Ab^p$ des $p$-relations ab\éliennes poss\ède une structure naturelle  d'espace vectoriel,  et l'application $(\omega_i)_i\mapsto (d\omega_i)_i$ permet  de d\éfinir 
une structure de  complexe sur ${Ab}^*=\oplus_{p\geq 0}Ab^p$  ($Ab^p=0$ pour $p>q$). Il est possible de donner des d\éfinitions globales, valables pour des tissus non n\écessairement compl\ètement d\écomposables, mais 
nous limiterons notre \étude au niveau des germes de tissu  en un point $m$ de $U$. 
  
 L'espace vectoriel  $Ab^p$  peut \^etre de dimension $r_p$ infinie. Mais, si $n$ est  un multiple de $q$,  A. H\énaut a d\émontr\é dans [H2] que cette dimension \était finie, et il en   a donn\é  une borne sup\érieure, g\én\éralisant le nombre  de Castelnuovo (r\ésultat d'abord  d\émontr\é par  Chern ([C]) si $p=q=1$, puis par  Chern-Griffiths ([CG]) pour $p=q$).  Ces auteurs, ainsi que Akivis ([A]), se sont particuli\èrement  int\éress\és \à la lin\éarisabilit\é des tissus, ainsi qu'\à la recherche de tissus de rang maximum, tels en particulier ceux  de $q$-rang maximum  propos\és par Goldberg ([G]).

 Lorsque $n$ n'est plus n\écessairement un multiple de $q$, Damiano ([D]) a  donn\é  une borne sup\érieure de $r_{n-1}$ pour les tissus en courbes ($q=n-1$). Il a en particulier montr\é  que le $(n+3)$-tissu en courbes qui g\én\éralise naturellement le tissu de Bol ([Bo]) du cas $n=2$, n'est pas  lin\éarisable, et est de rang $r_{n-1}$ maximum (pour $p<n-1$, les rangs $r_p$ de celui-ci sont  tous infinis). 
 
 Nous allons nous int\éresser \à des questions un peu diff\érentes. Toujours lorsque $n$ n'est plus n\écessairement un multiple de $q$, mais \à la condition que le tissu soit \emph{$p$-ordinaire}, nous allons montrer que $r_p$ est encore un nombre  fini, et nous pr\éciserons  sa  borne sup\érieure. Nous donnerons aussi une borne sup\érieure de la dimension $\widetilde r_p$ de l'espace des germes de $p$-relations ab\éliennes ferm\ées  pour les tissus \emph{fortement $p$-ordinaires} (la distinction n'\étant  \évidemment plus \à faire si $p=q$, et en particulier si $q=1$). 
 
 Nous allons pour cel\à 
 g\én\éraliser   aux tissus de   codimension  arbitraire les r\ésultats d\émontr\és dans  [CL] pour les tissus de codimension un. 
 Nous y avions   vu qu'en codimension un, le rang des $d$-tissus holomorphes  dans une vari\ét\é holomorphe de dimension $n$ \était major\é par un certain entier $\pi'(n,d)$,   strictement inf\érieur au nombre $\pi(n,d)$ de Castelnuovo
 pour $n\geq 3$, pourvu qu'ils soient \emph{ordinaires}.  Cette condition est  g\én\ériquement\footnote{Si $n=2$, tous les tissus sont ordinaires,  calibr\és, et $\pi'(2,d)=\pi(2,d)$.} v\érifi\ée. L'entier $\pi'(n,d)$ est donc aussi la borne sup\érieure du genre arithm\étique des courbes alg\ébriques \emph{ordinaires}\footnote{Ce sont les courbes dont les points d'intersection avec un hyperplan g\én\érique sont en ``position g\én\érale'', ou -de fa\c con \équivalente- dont le tissu associ\é dans l'espace projectif dual est ordinaire.}, et  accessoirement la borne inf\érieure du genre arithm\étique des courbes alg\ébriques arithm\étiquement de Cohen-Macaulay (cf. [GHL]). En outre, lorsqu'il existe un entier $k_0\geq 1$ tel que  $d$ soit  \égal \à la dimension $c(n,k_0)$ de l'espace des polyn\^omes homog\ènes de degr\é $k_0$ en $n$ variables (on dit alors que le tissu est \emph{calibr\é$^{1}$}), nous avions d\éfini -pour ces tissus- une  connexion holomorphe sur un certain fibr\é holomorphe de rang $\pi'(n,d)$, dont la courbure \était l'obstruction \à ce que le rang du tissu soit maximal \égal \à $\pi'(n,d)$, g\én\éralisant ainsi la courbure\footnote{Dans [DL1], nous avons aussi  propos\é un  programme sur Maple  pour calculer cette courbure, $n$ et $d$ \étant arbitraires, (une programmation avait d\éj\à \ét\é propos\ée par Pirio dans le cas $n=2$).} 
d\éfinie, pour $n=2$, d'abord par Pantazi ([Pa]), puis repris ind\épendamment par H\énaut ([H1]) et Pirio ([P]) en termes de connexions (c'est la courbure de Blaschke ([B]) si  $d=3$). 

Une fois surmont\ées   quelques difficult\és techniques suppl\émentaires, la m\éthode   est  la m\^eme en codimension arbitraire. Nous la r\ésumons ci-dessous pour les relations ab\éliennes, mais on proc\ède de la  m\^eme fa\c con pour les relations ab\éliennes ferm\ées, ainsi que nous le verrons dans la section 4 :

-  on observe d'abord  que les $p$-relations ab\éliennes  formelles \à un certain ordre $k$,  qui se projettent sur une relation ab\élienne  formelle $a_{k-1}$  donn\ée \à l'ordre $k-1$,  sont solution d'un syst\ème lin\éaire $\Sigma_k(a_{k-1})$ (sans second membre si $k=0$, et dont la partie homog\ène ne d\épend pas de $a_{k-1}$ si $k\geq 1$);   on \évalue la taille de ce syst\ème :

 - dire que le tissu est \emph{$p$-ordinaire} signifie que tous ces syst\èmes sont de rang maximum ; cette propri\ét\é ne requiert en fait qu'un nombre fini de conditions : plus pr\écis\ément, il suffit que les syst\èmes $\Sigma_k(a_{k-1})$ soient de rang maximum pour $k$ inf\érieur ou \égal \à un certain entier (not\é $k_p^0+1$ en g\én\éral, et  $k_p^0$ si $\Sigma_{k_p^0}$ a autant d'inconnues que d'\équations),  pour que le tissu soit $p$-ordinaire ;

- l'entier $k_p^0$ est  caract\éris\é par le fait   que les syst\èmes $\Sigma_k(a_{k-1})$ sont ou non  sur-d\étermin\és selon que $k>k^0_p$ ou $k\leq k^0_p$ ; on en d\éduit que  les   $p$-relations ab\éliennes   formelles \à l'ordre $k$  constituent, si  $k\leq k^0_p$, un fibr\é vectoriel $ R_k^p$ dont   le rang augmente avec $k$ jusqu'\à une certaine valeur   $\pi^0_p(n,d,q)$ que l'on sait calculer ;
 
 - puique  les syst\èmes lin\éaires pr\éc\édents sont sur-d\étermin\és pour $k>k_p^0$, la dimension de l'espace vectoriel des $p$-relations ab\éliennes formelles d'ordre $\infty$ en un point est au plus \égale \à $\pi^0_p(n,d,q)$ ; puisque le contexte est analytique, il en est  a fortiori de m\^eme pour le \emph{$p$-rang} du tissu (on appelle ainsi  la dimension maximum de l'espace des germes de $p$-relations ab\éliennes  en un point).

 Dire que le tissu est \emph{$p$-calibr\é} signifie que  le nombre d'\équations dans $\Sigma_{k_p^0}$ est \égal au nombre d'inconnues. Si le tissu est de plus  $p$-ordinaire, la projection $ R^p_{k_p^0}\to  R^p_{k_p^0-1}$  est alors  un isomorphisme. L'isomorphisme inverse permet de d\éfinir de fa\c con naturelle une connexion holomorphe sur le fibr\é  \break ${\cal E}:= R^p_{k_p^0-1}$, pour laquelle les sections \à d\ériv\ée covariante nulle s'identifient aux $p$-relations   
  ab\éliennes (m\éthode initi\ée dans [H1] quand  $n=2$, et utilis\ée dans [CL] en codimension un  pour $n$ quelconque). La courbure de cette connexion est donc une obstruction \à ce que le $p$-rang du tissu soit maximal \égal \à $\pi^0_p(n,d,q)$.  
Ainsi que nous l'avons fait dans [DL2] dans le cas de la codimension un, on pourrait -du moins en th\éorie- raffiner cette m\éthode des connexions   afin de  calculer explicitement   le $p$-rang du tissu,  y compris dans le cas non-calibr\é, sans avoir \à exhiber  les  relations ab\éliennes ; en pratique, il nous faudrait un ordinateur plus puissant pour arriver au bout  des calculs. 

On \étudie  de m\^eme le fibr\é $\widetilde  R_k^p$ des $p$-relations ab\éliennes  ferm\ées  formelles \à l'ordre $k$.

Observons par ailleurs  que les rangs $\rho_k(p)$ et $\widetilde \rho_k(p)$ des fibr\és $ R_k^p$ et $\widetilde  R_k^p$ (ainsi que la courbure dans le cas calibr\é) sont aussi des invariants des tissus, qui  peuvent suffire \à faire la distinction entre plusieurs d'entre eux, 
et qui sont  parfois  plus faciles \à calculer que les rangs  $r_p$ et $\widetilde r_p$ proprement dits. 

Dans la section 7, ind\épendante de ce qui pr\éc\ède,   nous montrons que l'un des exemples de  2-relation ab\élienne non-triviale exhib\é par Goldberg ([G]) pour un certain tissu de codimension 2 dans un espace de dimension 4,  n'est le  cobord d'aucune 1-relation ab\élienne, y compris au niveau des germes ; autrement dit, except\é pour $p=1$, il n'y a pas d'analogue au ``lemme de Poincar\é'' pour les relations ab\éliennes.

Je remercie vivement Alain H\énaut pour ses suggestions et encouragements.  
 
\section{D\éfinition des $p$-relations ab\éliennes \hb Notations et calculs pr\éliminaires  }
  
   \indent  Soient  
   
    $U$ un ouvert d'une vari\ét\é holomorphe de dimension $n$, ($n\geq 2$),
    
    $d$ un entier  $>0$, 
    
  et  $p$ et $q$ deux  entiers tels que   $ 0\leq p\leq  q\leq n-1$.

\n Notons   $T  {\cal F}  $  le fibr\é vectoriel holomorphe de rang $n-q$ des  vecteurs  tangents \à  un feuilletage holomorphe $ {\cal F} $ de codimension $q$.  Rappelons qu'une  $p$-forme   diff\érentielle  $\varpi$   sur $U$   est dite
 
   -   $ {\cal F} $-semi-basique si $i_X\varpi=0$ pour toute section  $X$ de  $T {\cal F} $, ($i_X$ d\ésignant le     produit int\érieur),  
   
   - $ {\cal F} $-invariante   si $L_X\varpi=0$ pour toute section  $X$ de  $T {\cal F}$, ($L_X=i_X\scirc d+d\scirc i_X$ d\ésignant la    d\ériv\ée de Lie), 
   
    - $ {\cal F} $-basique  si elle est \à la fois  $ {\cal F} $-semi-basique et $ {\cal F} $-invariante. 
    
    
    \n On notera  $B ^*({\cal F})$
    la sous-alg\èbre  diff\érentielle gradu\ée des formes  $ {\cal F}$-basiques.
    

    
     

\n Si $ {\cal F} $ est localement d\éfini par une submersion $\theta:U\to {\cal T}$ sur une vari\ét\é  ${\cal T}$ de dimension $q$, les formes  $ {\cal F} $-basiques sont aussi d\éfinies   localement comme \étant les images r\éciproques par $\theta$ des formes sur ${\cal T}$ ;  cette deuxi\ème d\éfinition, qui  ne d\épend pas de la submersion $\theta$ utilis\ée, permet de voir que le lemme de Poincar\é (localement, toute forme ferm\ée est exacte) est valable pour les formes $ {\cal F} $-basiques.
   
 On appellera \emph{syst\ème g\én\érateur de fonctions $ {\cal F} $-basiques} toute famille $u=(u_1,\cdots,u_q)$ de $q$ fonctions holomorphes ${\cal F}$-invariantes, telle que les fonctions $u_\alpha$   $(1\leq \alpha\leq q)$  et les 1-formes holomorphes $du_\alpha$ engendrent toute l'alg\èbre gradu\ée $B ^*({\cal F})$ : 
$du_1\wedge du_2\wedge\cdots\wedge du_q\neq 0$.

 On se donne  un $d$-tissu    holomorphe $\cal W$ de codimension $q$ sur $U$. 
     Bien qu'il soit possible de donner des d\éfinitions globales, on supposera $U$   suffisamment petit pour permettre des calculs locaux. En particulier, on supposera  
      que  $\cal W$ est compl\ètement d\écomposable  en  une famille $({\cal F}_i)_{1\leq i\leq d}$ de $d$ feuilletages holomorphes $({\cal F}_i)_{1\leq i\leq d}$  sur  $U$, en ``position g\én\érale".
     
     

\n {\bf D\éfinition 1 :}  

{\it   Soit $p$ un entier compris entre  0 et $q$. 
 
 Une  $p$-\emph{relation ab\élienne}  sur un  $d$-tissu $({\cal F}_i)_i$ de  codimension $q$ sur  $U$ est la donn\ée d'une famille de 
$p$-formes $(\omega_i)_i$  sur $U$, $(1\leq i\leq d)$, 
   
   - v\érifiant $\sum_i \omega_i\equiv 0 $ $($condition dite ``de trace nulle''$)$  ,
  
   - et telle que chaque forme $\omega_i$ soit ${\cal F}_i$-basique.

 \n   Si l'on impose en plus aux    formes ${\cal F}_i$-basiques $\omega_i$ d'\^etre ferm\ées, on dira que la $p$-relation ab\élienne  $(\omega_i)_i$ est   \emph{ferm\ée}.}

\n {\bf Remarques :} 

$(i)$ Pour $q=1$, les 1-relations ab\éliennes sont les relations ab\éliennes usuelles. La condition   de fermeture des formes $\omega_i$ est  alors automatiquement v\érifi\ée, comme c'est encore  plus g\én\éralement le cas si $p=q$. 
   
$(ii)$ L'ensemble des relations ab\éliennes  sur $U$ (resp. des  $p$-relations  ab\éliennes ferm\ées) poss\éde une structure naturelle d'espace vectoriel gradu\é.     $Ab^*(U)$ $\bigl($resp. ${\widetilde Ab ^*}(U)\bigr)$, et
l'on peut donner  des d\éfinitions analogues pour les germes de tissus en un point $m\in U$ et d\éfinir les espaces vectoriels gradu\és  $Ab^*_m$ et $\widetilde Ab^*_m$ des germes de  $p$-relations ab\éliennes \éventuellement ferm\ées.

 $(iii)$   L'application  $(\omega_i)_i\to  (d\omega_i)_i$  permet de d\éfinir une structure de  complexe sur les espaces  vectoriels  gradu\és $Ab^*(U)$ et $Ab^*_m$, dont on notera   
$H^*_{Ab}({\cal W})$ la cohomologie (sur un ouvert, ou au niveau des germes selon le contexte).

$(iv)$ 
Il est clair, si l'ouvert est connexe ou au niveau des germes,  que  $$H_{Ab}^0({{\cal W}})\cong\hbox{\bb C}^{d-1},$$
puisque c'est le noyau de l'application  $(k_1,k_2,\cdots, ,k_d)\to \sum_ik_i$ de 
{\bb C}$^d$ dans {\bb C}. 

$(v)$  Il est non moins clair, si l'ouvert est simplement  connexe, ou au niveau des germes,  que $$H_{Ab}^1({\cal W})=0,$$ puisque toute 1-relation ab\élienne ferm\ée $(\omega_i)_i$ se rel\ève par $d$ en une famille $(u_i)_i$ de fonctions basiques, en vertu du lemme de Poincar\é qui s'applique aux formes basiques de chaque feuilletage,  que $\sum_iu_i$ est une constante puisque $\sum_i \omega_i=0$, et que l'on peut toujours supposer cette constante nulle puisque les $u_i$ ne sont d\éfinies qu'\à une constante additive pr\ès. Cependant, 
le raisonnement pr\éc\édent ne se g\én\éralise pas  \à $H_{Ab}^p({\cal W})$ si $p>1$ :  nous  verrons dans la section 7 un exemple pour lequel  $H_{Ab}^2({\cal W})\neq 0$.
 
La  dimension, \éventuellement  infinie, de  $Ab^*_m\ \bigl($ resp. de $\widetilde Ab^*_m\bigr)$ s'appelle  {\it le $p$-rang} du tissu en $m$ (resp. {\it le $p$-rang ferm\é}),
   et sera not\ée $r_p(m)$ $\bigl($resp. $\widetilde r_p(m)\bigr)$. On appellera \emph{$p$-rang du tissu} (resp. \emph{$p$-rang fort}) la borne sup\érieure de ces nombres quand $m$ parcourt $U$.

 \n {\bf Proposition 1 :} {\it Si, parmi les $d$ feuilletages du tissu, il en existe deux, ${\cal F}_i$ et ${\cal F}_j$ $(i\neq j)$,  contenus dans un m\^eme  feuilletages ${\cal G}$ de codimension $q'$, $(1\leq q'<q)$, les p-rangs $r_p$ et $\widetilde r_p$ du tissu sont infinis  pour tout $p\leq q'$.}\hb
  \n {\it D\émonstration : } Toute $p$-forme  $\cal G$-basique $\varpi$ (\éventuellement ferm\ée) est en effet \à la fois  ${\cal F}_i$-basique et  ${\cal F}_j$-basique.   On  obtient par cons\équent   une p-relation ab\élienne   du tissu $(\omega_i)_i$ (\éventuellement ferm\ée) avec :
 $$\omega_i=\varpi\ ,\ \ \omega_j=-\varpi ,\ \ \omega_k=0  \hbox{ pour }k\neq i,j.$$
 
 \n {\bf Remarque  : } Si $d\varpi=0$, la classe de cohomologie des relations ab\éliennes ferm\ées   ainsi d\éfinies est nulle, puisqu'il existe une forme $\cal G$-basique $\eta$ telle que $d\eta=\varpi$, et que la relation ab\élienne pr\éc\édente  est la diff\érentielle de 
 $(\eta_i=\eta\ ,\ \ \eta_j=-\eta ,\ \ \eta_k=0  \hbox{ pour }k\neq i,j)$.
 
 \rightline{QED} 
 
 Pour \éviter d'utiliser  trop d'espace, on notera d\ésormais ${\bf b(r,s)}$  le coefficient binomial $\frac{r!}{s!(r-s)!} $,
au lieu de la notation usuelle $\begin{pmatrix}r\\s\end{pmatrix} $. 
  On notera aussi 
  {\bf c(r,h)}$:=b(r-1+h,h)$ la dimension de l'espace vectoriel des polyn\^omes  homog\ènes de degr\é $h$, \à $r$ ind\étermin\ées,  \à coefficients dans {\bb C}.
 
 Pour tout entier $p$  $ (1\leq p\leq q)$, notons  
  ${\cal A}_p $ (resp. ${\cal B}_p $) l'ensemble des multi-indices 
 $$A=(1\leq \alpha_1< \alpha_2<\cdots< \alpha_p\leq q),\  
  \bigl(\hbox{resp.  } B=(1\leq \lambda_1< \lambda_2<\cdots< \lambda_p\leq n)\bigr),$$

Soit  $S(r,h)$   l'ensemble des $c(r,h) $ multi-indices $L$ de d\érivation d'ordre $h$ des fonctions holomorphes de $r$ variables
 $$L=(\rho_1,\cdots,\rho_r)\ , \hskip 1cm \rho_j \geq 0\ ,\hskip 1cm \sum_{\rho=1}^r \rho_j =h \ (\hbox{encore not\é }|L| ),$$ et $v'_L$ la d\ériv\ée correspondante d'ordre $h$ d'une telle fonction (si $|L|=0$, on notera $o_r$ le multi-indice correspondant, et l'on conviendra que $v'_{o_r}=v$). 
 
 Si $L=(\rho_1,\cdots,\rho_r)$ et $L'=(\rho'_1,\cdots,\rho'_r)$, $L+L'$ d\ésignera le multi-indice 
 $(\rho_1+\rho'_1,\cdots,\rho_r+\rho'_r)$. 
 
 Notons $L+1_j$ le multi-indice obtenu \à partir de $L$ en augmentant $\rho_j$ d'une unit\é ; si $\rho_j\geq 1$, on notera aussi $L-1_j$ le multi-indice obtenu \à partir de $L$ en diminuant $\rho_j$ d'une unit\é (si $\rho_j=0$, on conviendra que $v'_{L-1_j}=0$). On notera $o_r$ le multi-indice $(\rho_1=0,\cdots,\rho_r=0)$ et l'on posera $v'_{o_r}=v$.

\n {\bf Lemme 1  :} {\it   Soient

  $\cal F$ un un feuilletage holomorphe de codimension $q$ sur $U$, 

 $u=(u_1,\cdots, u_q)$ un syst\éme g\én\érateur de $q$ fonctions holomorphes sur $U$, $\cal F$-invariantes, 
 
 $f$ une fonction holomorphe de $q$ variables, et $J$ une fonction holomorphe sur $U$.
 
\n La formule suivante est alors v\érifi\ée pour tout  $L\in S(n,k)$, $k\geq 1$ :
$$\bigl((f\scirc u).J\bigr)'_L=\sum_{h=0}^k\ \ \sum_{K\in S(q,h)} M_{L}^{K}({\cal F},J) .\bigl((f)'_K\scirc u\bigr),$$
 les coefficients  $  M_{L}^{K}({\cal F},J)$, not\és  $M_{L}^{K}$ s'il n'y a pas d'ambigu\"\i t\é et qui ne d\épendent pas de $f$, 
  \étant des fonctions holomorphes d\éfinies par r\écurrence sur $k$ par les formules suivantes :
  
  $M^{o_q}_{1_\lambda}=J'_\lambda$  \hskip 1cm et \hskip 1cm   $M^{1_\alpha}_{1_\lambda}=(u_\alpha)'_\lambda .J$ \hskip 1cm pour $k=1$,
 
 
 
 \n et pour tout multi-indice $L=(\ell_1,\cdots,\ell_n)$    de degr\é $|L|=k\geq 1$ : 
 
 $ M_{L}^{o_q}=J'_L$

$ M_{L+1_\lambda}^{K}=( M_{L}^{K})'_\lambda+\sum_{\alpha} M_{L}^{K-1_\alpha}.(u_\alpha)'_\lambda$ \hskip 2cm pour $1\leq  |K|\leq k$,

$M_{L+1_\lambda}^{K}=\sum_{\alpha} M_{L}^{K-1_\alpha}.(u_\alpha)'_\lambda$ \hskip 2cm pour $|K|=k+1$,

 \n $($Il est sous-entendu, si $K=(k_1,\cdots,k_\alpha,\cdots,k_q)$,  que toutes les sommations $\sum_\alpha$ ci-dessus  sont limit\ées aux couples $(\alpha,K)$ tels que $ k_\alpha>0).$
 
 \n Pour $J\equiv 1$  on posera  : $N_{L}^{K}({\cal F}) :=M_{L}^{K}({\cal F},1)$.

 En particulier, pour  $|K|=|L]$, notant  $s$ la valeur commune de ces deux entiers, et posant :\hb  $L=(\ell_1,\cdots,\ell_\lambda,\cdots,\ell_n)$, $K=(k_1,\cdots,k_\alpha,\cdots,k_q)$, on obtient  :
$$N_{L}^{K}({\cal F})=\sum\Bigl( (u_{\alpha_1})'_{\lambda_1}.(u_{\alpha_2})'_{\lambda_2}\cdots(u_{\alpha_s})'_{\lambda_s}\Bigr),\hbox{ et }M_{L}^{K}({\cal F},J)=J.N_{L}^{K}({\cal F}), $$  la sommation \étant effectu\ée sur tous les produits
$\prod(u_\alpha)'_\lambda $ de $s$ d\ériv\ées premi\ères $(u_\alpha)'_\lambda$ tels que  
 l'entier $\alpha$ $(1\leq \alpha\leq q)$ figure $k_\alpha$ fois dans la suite $(\alpha_1,\cdots,\alpha_s)$, et l'entier $\lambda$ $(1\leq \lambda\leq n)$ figure $\ell_\lambda$ fois dans la suite  $(\lambda_1,\cdots,\lambda_s)$}.\hb 
 {\it D\émonstration :} Ces formules s'obtiennent par r\écurrence sur $|L|$, en d\érivant $\bigl((f\scirc u).J\bigr)'_L$ par rapport \à $x_\lambda$ :
$$\Bigl(\sum_{h=0}^k\ \ \sum_{K\in S(q,h)}M_{L}^{K} .\bigl(f'_K\scirc u\bigr)\Bigr)'_\lambda=
\sum_{h=0}^k\ \ \sum_{K\in S(q,h)}\Bigl[( M_{L}^{K})'_\lambda .\bigl(f'_K\scirc u\bigr)+ M_{L}^{K}.\sum_{\alpha=1}^q
(u_\alpha)'_\lambda.\bigl(f'_{K+1_\alpha}\scirc u\bigr)\Bigr]
.$$
\n  {\bf Rappels sur les complexes de Spencer et de Koszul }

Notant $T^*U$ le fibr\é tangent complexe de $U$, le complexe de de Rham $(\Omega^*(U),d)$ des formes holomorphes sur $U$ induit, pour tout entier positif $k$, le  complexe  de Spencer $(Sp_k)$, et le  noyau de la projection $(Sp_k)\to(Sp_{k-1})$ n'est autre que le complexe de Koszul $(Ko_k)$, toujours acyclique :
$$ \begin{matrix}
&&&0&&0&&0 \\
&&&\downarrow&&\downarrow&&\downarrow \\
(Ko_k)&&\cdots\to&S^{k+1} T^*(U)\otimes \wedge^{p-1} T^*(U)&\buildrel{d_{p-1}}\over\longrightarrow&  S^{k} T^*(U)\otimes \wedge^{p} T^*(U)&\buildrel{d_p}\over\longrightarrow&S^{k-1} T^*(U)\otimes \wedge^{p+1} T^*(U) \\
\downarrow&&&\downarrow&&\downarrow&&\downarrow \\
(Sp_k)&\hskip .2cm  &\cdots\to &J^{k+1}\Bigl(\bigwedge^{p-1}T^*U\Bigr)&\to & J^{k}\Bigl(\bigwedge^{p}T^*U\Bigr)&\to&  J^{k-1}\Bigl(\bigwedge^{p+1}T^*U\Bigr) 
  \\
\downarrow&&&\downarrow&&\downarrow&&\downarrow \\
 
(Sp_{k-1})&\hskip .2cm  &\cdots\to &J^{k}\Bigl(\bigwedge^{p-1}T^*U\Bigr)&\to & J^{k-1}\Bigl(\bigwedge^{p}T^*U\Bigr)&\to&  J^{k-2}\Bigl(\bigwedge^{p+1}T^*U\Bigr) 
  \end{matrix} $$



\n {\bf Lemme 2 :} {\it 

\n $(i)$ Pour tout entier $p$, $(1\leq p\leq n)$, les $k$-jets de $p$-formes ferm\ées en un point $m$ de $U$ se projetant sur un $(k-1)$-jet donn\é forment un espace affine  de dimension 
$$ {\bf z(n,p,k)}:=\sum_{j=p}^{n}(-1)^{p-j}\ b(n,j).c(n,k+p-j),\ 
\Bigl( =\sum_{j=0}^{p-1}(-1)^{p-j+1}\ b(n,j).c(n,k+p-j)\Bigr) .$$

\n $(ii)$ }  l'identit\é suivante est v\érifi\ée :
  $$ z(n,p,k)\equiv \ b(n+k,n-p)\ .\ c(p,k)$$
 
\n {\it D\émonstration :} 

Une fois fix\ée la famille $\Bigl((g_B)'_L. dx_B)\Bigr)_{L\in S(n,h),h\leq k-1}$ d\éfinissant le $k-1$-jet, le $k$-jet de la forme $\omega=\sum_{B\in {\cal B}_p} g_B. dx_B$ est d\éfini par la famille $$\Bigl((g_B)'_L. dx_B)\Bigr)_{L\in S(n,k)}.$$ Si $d\omega=0$, cet \él\ément de $ S^k T^*(U)\otimes \wedge^p T^*(U)$ appartient au noyau de  l'application $$ d_p:S^k T^*(U)\otimes \wedge^p T^*(U)\rightarrow S^{k-1} T^*(U)\otimes \wedge^{p+1} T^*(U)$$ du complexe de Koszul, 
et r\éciproquement. Par cons\équent la dimension de l'espace de ces \él\éments est \égale \à la somme altern\ée des dimensions des termes \à partir de l\à dans le complexe de Koszul,  ceux-ci formant une r\ésolution de ce noyau.  Il revient au m\^eme,  puisque le complexe de Koszul  est acyclique, de prendre la somme altern\ée des dimensions des termes   pr\éc\édant $ S^k T^*(U)\otimes \wedge^p T^*(U)$) dans ce complexe. On en d\éduit la partie $(i)$ du lemme.

L'expression  $b(n+k,n-p)\ .\ c(p,k)$ est  \égale \à $$\frac{1}{(p-1)!(n-p)!}\times\prod_{1\leq s\leq n, \ s\neq p}(k+s). $$Elle est en particulier  polyn\^omiale de degr\é $n-1$ en $k$. Or 
chaque expression $c(n,k+p-j)$ est \également polynomiale de degr\é $n-1$  en $k$, donc aussi $z(n,p,k)$ que l'on peut ainsi prolonger  aux valeurs n\égatives de $k$.  Puisque 
$$c(n,k+p-j)=\frac{1}{(n-1)!}\prod_{1\leq r\leq n-1}(k-p-j+r)\ ,$$
$c(n\ , -s+p-j)=0$ pour $0\leq j\leq p-1$ si $s$ est un entier compris entre $p+ 1$ et $n$, et  pour $p\leq j\leq n$ si $s$ est un entier compris entre $1$ et $p-1$. 
Ceci prouve  que $z(n,p\ ,-s)=0$ pour tout entier $s\in [1,n]$, \ $s\neq p$. Les deux polyn\^omes en $k$ de degr\é $n-1$,  $z(n,p,k)$ et $b(n+k,n-p).c(p,k)$, ont donc les m\^emes racines. Ils ont aussi le m\^eme terme constant $b(n,p)$. Ils sont donc \égaux, d'o\ù $(ii)$.

\n {\bf Corollaire :} {\it Les $k$-jets de $p$-formes ferm\ées $\cal F$-basiques en un point $m$ de $U$,  se projetant sur un $(k-1)$-jet donn\é $(\cal F$ d\ésignant un feuilletage de codimension $q)$, forment un espace affine de dimension $z(q,p,k)$.
 }
\n {\it D\émonstration :} Ces formes s'identifient  en effet naturellement  aux $p$-formes ferm\ées sur une sous-vari\ét\é  de dimension $q$ transverse \à $\cal F$. 

\section{Relations ab\éliennes des tissus ordinaires}

 Rappelons le r\ésultat suivant :
 
 \n{\bf Th\éor\éme 1 } {\it  $\bigl($A.  H\énaut $( [$H2$])\bigr)$ :} 
   
  {\it Supposons 
que  les feuilles du tissu sont en position g\én\érale en tout point $m$ de $U$,  
  et supposons  de plus que    $n$ est un multiple de $q$. Alors  :
  
$(i)$  L'espace vectoriel    $Ab_m^p$ des germes de $p$-relations ab\éliennes  en un  point $m$ de $U$ a une  dimension finie  $r_p$ qui ne  d\épend pas de  $m $. 

$(ii)$ Cette dimension   $r_p$   est major\ée par le nombre
$$ \pi_p(n,d,q) = b(q,p).    
  \sum_{h\geq 0} c(q,h). \biggl(d-\Bigl(\frac{n}{q}-1\Bigr)(p+h) -1\biggr)^+  ,  
  $$
  \indent la  notation $a^+$ d\ésignant,  pour tout nombre r\éel  $a$, le nombre    sup $(a,0)$.
  
$(iii)$ Cette borne est optimale : il existe un $d$-tissu sur $U$, dont le rang   $r_p$ est \égal \à $\pi_p(n,d,q)$.}
 
\n {\bf Remarque :} Ces  nombres  $\pi_p(n,d,q)$ g\én\éralisent les nombres   $\pi(n,d)=\pi_1(n,d,1)$ de Castelnuovo. Et ce  th\éor\éme  \était  d\éj\à d\émontr\é  par Chern ([C]) pour    $q=1$, et plus g\én\éralement par Chern-Griffiths   ([CG]) lorsque  $p=q$.

Revenons au cas g\én\éral o\ù $n$ n'est plus n\écessairement un multiple de $q$. 
On se donne un $d$-tissu    de codimension $q$ sur $U$, que l'on suppose compl\étement d\écomposable et d\éfini par $d$ feuilletages $({\cal F}_i)_{1\leq i\leq d}$.  Soit $(x_1,\cdots,x_\lambda,\cdots,x_n)$ un syst\éme de coordonn\ées locales sur $U$, et pour tout $i=1,\cdots,d$,  on note $u_i:=(u_{i,1},\cdots,u_{i,q})$ un syst\éme g\én\érateur de fonctions ${\cal F}_i$-invariantes.

 \n   Toute $p$-forme ${\cal F}_i$-basique $$\omega_i=\sum_{A\in {\cal A}_p} (f_{i,A}\scirc u_i)\  du_{i,A} $$ 
 s'\écrit encore : 
 $$\omega_i=\sum_{B\in {\cal B}_p}\sum_{A\in {\cal A}_p} (f_{i,A}\scirc u_i)\ J_{i,B}^A\  dx_{B} ,$$
$J_{i,B}^A$ d\ésignant le d\éterminant de la matrice jacobienne 
$$J_{i,B}^A:=det\biggl(\frac{D(u_{i,\alpha_1},\cdots,u_{i,\alpha_p})}{D(x_{\lambda_1},\cdots,x_{\lambda_p})}\biggr).$$
 Ainsi, les $p$-relations ab\éliennes  s'identifient aux familles $F=\Bigl(f_{i,A}\Bigr)_{i, A}$ de $d\times b(q,p)$ fonctions holomorphes $f_{i,A}$ de $q$ variables, satisfaisant aux $b(n,p)$ identit\és
 $$(E_B)\hskip 2cm \sum_{i=1}^d\sum_{A\in {\cal A}_p} (f_{i,A}\scirc u_i)\ J_{i,B}^A\equiv 0. $$
 \n Puisque $\Bigl(\sum_{i=1}^d\sum_{A\in {\cal A}_p}( f_{i,A}\scirc u_i).J_{i,B}^A\Bigr)'_L=\sum_{i=1}^d\sum_{A\in {\cal A}_p}\sum_{ |K|\leq  |L|} M_{B,L}^{i,A,K}.\Bigl((f_{i,A})'_K\scirc u_i\Bigr),$\hb
 les identit\és $$(E_B)'_L
  \hskip 2cm \Bigl(\sum_{i=1}^d\sum_{A\in {\cal A}_p}( f_{i,A}\scirc u_i).J_{i,B}^A\Bigr)'_L\equiv 0,\hskip 1cm L\in S(n,h),\hskip .3cm 0\leq h\leq k ,$$s'\écrivent encore 
 $$\sum_{i=1}^d\sum_{A\in {\cal A}_p}\sum_{ |K|\leq  |L|} M_{B,L}^{i,A,K}.\Bigl((f_{i,A})'_K\scirc u_i\Bigr)\equiv 0, $$ 
 o\ù l'on a pos\é $$M_{B,L}^{i,A,K}:=M_{L}^{K}({\cal F}_i, J_{i,B}^A)$$. 
 
 \n  On ordonne  de 1 \à $b(q,p)$ (resp.  $b(n,p)$) les \él\éments   de  $ {\cal A}_p$ (resp. $ {\cal B}_p$).
  On ordonne  de m\^eme de 1 \à $c(q,k)$ (resp.  $c(n,k)$) les \él\éments   de  $ S(q,k)$ (resp. $S(n,k)$).
  On ordonne alors les indices $(A,K)$ (resp. $(B,L)$ suivant l'ordre lexicographique :
  
  $(A,K)<(A',K')$,    si $A<A'$ ou si ($A=A'$ et $K<K'$), et  r\ègle analogue  pour les indices $(B,L)$.
  
 \n  Quant aux indices $(i,A,K)$ on les ordonne suivant la r\ègle :
  
 $(i,A,K)<(i',A',K')$ si  $(A,K)<(A',K')$, ou si $(A,K)=A',K')$ et $i<i'$.

 \n Notons :
  
  $\Theta ^r$ le fibr\é vectoriel holomorphe trivial de rank $r$,
  
  $  \beta_k(p):=\sum_{h=0}^k b( n,p).c(n,h),$
  
   $ \alpha_k(p):=d. \sum_{h=0}^k  b(q,p). c(q,h),$
   
    $P^{(k)}_h(p)$ la matrice $\bigl(\!\bigl(M_{(B,L)}^{(i,A,K)}\bigr)\!\bigr) $, de taille $b( n,p).c(n,k)\times d. b(q,p). c(q,h)$,  obtenue  pour $|L|=k$ et  $|K|=h$ (avec la convention $P^{(k)}_h=0$ si $h>|L|$).  On \écrira aussi $P _k(p)$ (voire $P_k$ s'il n'y a pas d'ambigu\"\i t\é sur $p$),  au lieu de $P^{(k)}_k(p)$.
    
     ${\cal M}_k(p)$ la matrice de taille $ \beta_k(p)\times  \alpha_k(p)$ construite avec les blocs ${P}_h^{(\ell)} $
 ($h,\ell\leq k$), le  bloc ${P}_{h+1}^{(\ell)}$ \étant \à droite de  ${P}_{h}^{(\ell)}$, et ${P}_{h+1}^{(\ell)}$ en dessous, 
 
et $Q_{k}(p)$ la sous-matrice  de taille   $\bigl(b( n,p).c(n,k)\bigr)\times \alpha_{k-1}(p)$  dans  ${\cal M}_{k}(p)$ form\ée avec les blocs  ${P}_{h}^{(k)}$ pour     $0\leq h\leq k-1$ :
 
$$\begin{matrix}{\cal M}_k(p)   &  =   &
 \begin{pmatrix}{P}_0^{(0)}=P_0&0&0&....&....&0&0\\ 
 {P}_0^{(1)}& {P}_1^{(1)}=P_1&0&....&....&0&0&\\
....& .... &....&....&....&0&0&\\ 
  P^{(k-1)}_0& P^{(k-1)}_1&P^{(k-1)}_2&....&....&P^{(k-1)}_{k-1}=P_{k-1}&0\\ 
P^{(k)}_0& P^{(k)}_1&P^{(k)}_2&....&....&P^{(k)}_{k-1}&P^{(k)}_{k}=P_{k}\\ 
 & \end{pmatrix}
 \\ 
 
 Q_{k} (p)&=&\hskip -2.3cm \begin{pmatrix} P^{(k)}_0&\hskip .8cm P^{(k)}_1&\hskip .7cm P^{(k)}_2&....&....&\hskip .8cm P^{(k)}_{k-1}\end{pmatrix}

\hskip .3cm
 \end{matrix}$$
[On omettra parfois  la parenth\èse $(p)$ si  aucune ambigu\"\i t\é n'est \à craindre].

  
\n   Notons :
  
  $R_k^p$ l'ensemble des relations ab\éliennes  formelles \à l'ordre $k$, que l'on identifie localement \à un sous-ensemble de l'espace total du fibr\é trivial  $\Theta^{\alpha_k(p)}$,
  
   $R_k^p\to R_{k-1}^p$ la projection naturelle,
   
    $w(i,A,K)$ $\bigl(A\in {\cal A}_p$, $K\in S(q,h)$,\ $0\leq h\leq k\bigr)$ 
le  nombre  candidat  \à repr\ésenter la valeur en un point $u_i(m)$ d'une fonction $(f_{i,A})'_K$, 

$w_h$ le $d\times b(q,p).c(n,h)$-vecteur colonne des  $w(i,A,K)$ pour $|K|=h$,

et $w^{(k)}$ le $\alpha_k(p)$-vecteur colonne $(w_0,w_1,\cdots w_k)$ des  $w(i,A,K)$ pour  $|K|\leq k$.

  \n  On d\éduit de ce qui pr\éc\ède le 

  \n {\bf Th\éor\ème 2 :}  {\it Supposons   $1\leq p<q$.

\n $(i)$  Localement, $R_k^p$ s'identifie au noyau de  ${\cal M}_{k}(p)$ $\bigl($inclus  dans le fibr\é trivial  $\Theta^{\alpha_k(p)} \bigr)$. Si la matrice  ${\cal M}_{k}(p)$ conserve un rang  constant en tout point $m\in U$, $R_k\to U$ est un fibr\é vectoriel holomorphe de rang
  $$\rho_k(p):= \alpha_k(p)-rang\ {\cal M}_{k}(p).$$

\n $(ii)$ Les \él\éments de $R_k^p$ se projetant sur un \él\ément $a_{k-1}\in R_{k-1}^p$ donn\é sont   solution du syst\éme lin\éaire suivant avec second membre  :
$$  \Sigma_k(a_{k-1})\hskip 1cm <P_k,w_k>\ =\ -<Q_k(p),a_{k-1} >$$ des 
 $b(n,p).c(n,k)$ \équations  $(E_B)'_L$  ($B\in {\cal B}_p$, $L\in S(n,k)$)  \hb
et   $d. b(q,p). c(q,k)$ inconnues $w(i,A,K)$ $\bigl(A\in {\cal A}_p$, $K\in S(q,k)\bigr)$.}

\n {\bf D\éfinition 2 : }{\it  Le tissu est dit \emph{$p$-ordinaire} si, pour tout entier $k\geq 0$, la matrice $P_k$ est de rang maximum {\rm inf} $\bigl(b(n,p).c(n,k)\ ,\ d. b(q,p). c(q,k)\bigr)$}.

 \n {\bf Remarques :} 

1- Le syst\ème $\Sigma_k(a_{k-1})$  ayant une signification intrins\èque, cette d\éfinition ne d\épend pas des syst\èmes g\én\érateurs $u_i$ utilis\és, pas plus que des coordonn\ées locales $(x_1,\cdots,x_n)$.

2- Si  $p=q$ (et en particulier si $q=1$),  la $1$-forme $\sum_i\omega_i$  est automatiquement ferm\ée  : les \équations \à l'ordre $k$ sont alors  en nombre $z(n,q,k)$   et non   $b(n,q)\times c(n,k)$. Ce cas rel\ève donc  des relations ab\éliennes ferm\ées  trait\ées dans la section suivante. 

\n {\bf Lemme 3 : }{\it Le rapport  $ b(n,p).c(n,k)/ b(q,p). c(q,k)$ est une fonction strictement croissante de $k$ pour $k\geq 0$. }

\n {\it D\émonstration :} Le rapport pr\éc\édent s'\écrit en effet sous la forme 
$$\frac{b(n,p)}{b(q,p)}\ .\ \prod_{j=1}^{n-q} \frac{(q-1+j+k) }{(q-1+j) } .$$

\rightline{QED}

\n Soit alors $ {\bf k^0_p} $ le plus grand entier $k$ tel que $ b(n,p).c(n,k)/ b(q,p). c(q,k)$ soit au plus \égal \à $d$. 

 \n On obtient  le 
 \n {\bf Th\éor\ème 3 :}
{\it Supposons $1\leq p<q$, le $d$-tissu de codimension $q$ ci-dessus $p$-\emph{ordinaire},  et $d> \frac{b(n,p)}{b(q,p)}$. 
 
 $(i)$ Pour $k\leq k^0_p$, l'ensemble $R_k^p$ des $p$-relations ab\éliennes  formelles  \à l'ordre $k$ poss\ède une structure naturelle de fibr\é vectoriel holomorphe de rang 
$$b(q,p) \sum_{h=0}^k c(q,h) \Bigl(d-\frac{b(n,p).c(n ,h)}{b(q,p).c(q ,h)}\Bigr).$$

 $(ii)$ Le $p$-rang  $r_p$ du tissu  est major\é par le nombre 
 $$\pi^0_p(n,d,q):=b(q,p) \sum_{h\geq 0} c(q ,h) \biggl(d-\frac{b(n,p).c(n ,h)}{b(q,p).c(q ,h)}\biggr)^+.$$}
 

\n {\it D\émonstration :} 

Si 
le $d$-tissu de codimension $q$ ci-dessus est  \emph{$p$-ordinaire}, et si $h\leq k^0_p$, l'espace des  $p$-relations ab\éliennes  formelles  \à l'ordre $h$ se projetant sur une $p$-relation ab\élienne   formelle  \à l'ordre $h-1$
donn\ée,  est un espace affine de dimension \égale \à la diff\érence $$  d. b(q,p). c(q,h)-b(n,p).c(n,h)$$ du nombre d'inconnues et du nombre d'\équations, d'o\ù la partie $(i)$ du th\éor\ème. 

Si 
le $d$-tissu  est  $p$-\emph{ordinaire}, et si $h> k^0_p$, l'espace des  $p$-relations ab\éliennes  formelles \à l'ordre $h$ se projetant sur une $p$-relation ab\élienne   formelle \à l'ordre $h-1$
donn\ée,  est un espace affine de dimension 0 ou est vide. L'espace des $p$-relations ab\éliennes  formelles    d'ordre $\infty$ est donc de rang au plus \égal \à celui $\pi^0_p(n,d,q)$ de $R^p_{k_0(p)}$. Il en est de m\^eme pour l'espace des germes de $p$-relations ab\éliennes    en un point de $U$, puisque le contexte est analytique, d'o\ù la partie $(ii)$ du th\éor\ème. 


\rightline{QED}

 \n {\bf Remarques :} 
 
 $(i)$ Lorsque  $n$ est un multiple de $q$, $\pi^0_p(n,d,q)$ est strictement plus petit que $\pi_p(n,d,q)$, sauf pour $n=2$, $q=1$ (auquel cas il y a \égalit\é).
 
 $(ii)$ Quand il existe un $d$-tissu $p$-ordinaire \emph{parall\élisable} (cf. section 6) de codimension $q$ dans un espace ambiant $n$-dimensionnel, celui-ci a un $p$-rang maximal, et cette borne est donc optimale. 

\n {\bf Th\éor\ème 4   : } {\it Pour que le tissu soit   $p$-ordinaire,  il suffit que  $ P_k(p)$ soit  de rang maximum pour   $k\leq  k^0_p+1$ $($et m\^eme seulement pour $k\leq  k^0_p$ lorsque $ P_{k_p^0}(p)$ est une matrice carr\ée$)$.}

\n {\it D\émonstration :} 

Supposons $k>k^0_p$, et $P_k$  de rang maximum :   son rang est donc \égal au nombre de ses colonnes (on dira pour abr\éger que $P_k$ est une matrice``injective").

Le syst\ème $\Sigma_k$ du th\éor\ème 2 ci-dessus provient de ce que, pour toute  relation ab\élienne faible $F$, 
$$  <P_k,\Bigl((f_{i,A})'_K\scirc u_i\Bigr)_{|K|=k}  >\ \equiv\ - <Q_k(p),j^{k-1}F >.$$En d\érivant successivement par rapport aux diff\érentes coordonn\ées $x_\lambda$ (notation abr\ég\ée : $\partial_\lambda$), on voit que, pour tout $\lambda=1,\cdots, n$,  
 $$<P_k,\Bigl(\sum_\alpha\bigl((f_{i,A})'_{K+1_\alpha}\scirc u_i)(u_{i,\alpha})'_\lambda\Bigr)_{|K|=k}  >$$
 ne d\épend que de $j^k F$, et 
 $j^{k+1} F$ est   d\éfini de fa\c con unique par les donn\ées de $j^k F$ et de la famille $\Bigl(\partial_\lambda\bigl((f_{i,A})'_K\scirc  u_i \bigr)\Bigr)_{(i,A,|K|=k,\lambda)}$. Il en r\ésulte que $w_{k+1}=\bigl(w(i,A,K')\bigr)_{(i,A,|K'|=k+1)}$ est aussi bien d\éfini par la famille $$ \Bigl(\sum_\alpha\bigl(w(i,A,K+1_\alpha).(u_{i,\alpha})'_\lambda\Bigr)_{(i,A,|K|=k,\lambda)}.$$ Ainsi, d\ès lors que $P_k$ est une matrice injective,  $P_{k+1}$ l'est aussi. Puisque   $P_{k^0_p+1}$  est injective par hypoth\èse (resp. $ P_{k_p^0}$ si cette matrice est carr\ée), il en est de m\^eme pour tout  $P_k$,  $k>k^0_p$.
 
 \rightline{QED}

\section{Relations ab\éliennes ferm\ées}

Reprenons les notations de la section pr\éc\édente.
Pour obtenir  des $p$-relations ab\éliennes ferm\ées, il faut maintenant, dans les \équations de la section pr\éc\édente,  

\n 1)  remplacer la base $\Bigl(w(i,A,K)\Bigr)_{{A\in {\cal A}_p, |K|=k}}$ de $  S^{k} T^*(W_i)\otimes \wedge^{p} T^*(W_i)$ par une base du noyau de 
$d_p$ dans  le complexe de Koszul relatif \à une  sous-vari\ét\é $W_i$ de dimension $q$ transverse \à ${\cal F}_i$, 

\n 2) r\éduire \à $z(n,p,k)$ le nombre $b(n,p).c(n,k)$ des \équations $(E_B)'_L$ ($B\in {\cal B}_p$, $L\in S(n,k)$) : ces derni\ères, en effet,  ne sont plus   toutes  lin\éaitrement ind\épendantes, puisque $d_p(\sum_i w(i,k)=0$ dans le complexe de Koszul relatif \à $U$. 


Utilisant maintenant, pour tout $i=1,\cdots ,d$, une base (avec $z(q,p,k)$ \él\éments) de l'espace  des $k$-jets de $p$-formes ${\cal F}_i$-basiques \emph{ferm\ées},  
notons   $\widetilde{\cal M}_k(p)$ la matrice de taille $ \tilde\beta_k(p)\times  \tilde\alpha_k(p)$ correspondante, construite de fa\c con analogue \à  ${\cal M}_k(p)$ dans la section pr\c\édente, avec des blocs $\widetilde {P}_h^{(\ell)} $ de taille \hb $z(n,p,\ell)\times z(q,p,h)$,
 o\ù 
   
    $\tilde \alpha_k(p):=d.\sum_{h=0}^k  z(q,p,h),$
    
    et  $  \tilde{\beta}_k(p):=\sum_{h=0}^k z(n,p,k).$

\n On d\éfinit de m\^eme $ \widetilde Q_{k}(p)$,  de taille   $z(n,k)\times \tilde\alpha_{k-1}(p)$  dans  $\widetilde{\cal M}_{k}(p)$ form\ée avec les blocs  $\widetilde{P}_{h}^{(k)}$ pour     $0\leq h\leq k-1$, et l'on \écrira aussi $\widetilde {P}_k$ au lieu de $\widetilde {P}_k^{(k)} $.

Notant maintenant  $\widetilde{R}_k^p$ l'ensemble des relations ab\éliennes ferm\ées   formelles \à l'ordre $k$ et\break $\widetilde{R}_k^p\to \widetilde{R}_{k-1}^p$ la projection naturelle, on obtient  le 
 
  \n {\bf Th\éor\ème 5 :}  {\it Si  $q\geq 1$,

\n $(i)$  Localement, $\widetilde R_k^p$ s'identifie au noyau de  $\widetilde{\cal M}_{k}(p)$ $\bigl($inclus  dans le fibr\é trivial  $\Theta^{\tilde\alpha_k(p)} \bigr)$. Si la matrice  $\widetilde{\cal M}_{k}(p)$ conserve un rang  constant en tout point $m\in U$, $\tilde R_k\to U$ est un fibr\é vectoriel holomorphe de rang
  $$\tilde\rho_k(p):=\tilde \alpha_k(p)-rang\ \widetilde{\cal M}_{k}(p).$$

\n $(ii)$ Les \él\éments de $\widetilde R_k^p$ se projetant sur un \él\ément $\widetilde a_{k-1}\in \widetilde R_{k-1}^p$ donn\é sont   solution d'un  syst\éme lin\éaire avec second membre $\widetilde \Sigma_k(\widetilde a_{k-1})$ \hb
 - de $z(n,p,k)$ \équations  $(E_B)'_L$    \hb
- et   $d. z(q,p,k)$ inconnues, 
\n que l'on peut \écrire en abr\ég\é :
$$<\widetilde P_k,\tilde w_k>\ =\ < \widetilde Q_k(p),\widetilde a_{k-1}>.$$ }
 {\bf D\éfinition 3 : }{\it  Le tissu est dit \emph{$p$-fortement ordinaire} si, pour tout entier $k\geq 0$, la matrice $\widetilde P_k$ est de rang maximum \emph{inf }$\bigl(z(n,p,k)\ ,\ d. z(q,p,k)\bigr)$}. 

\noindent Cette d\éfinition a une signification intrins\èque, comme la d\éfinition 2.

 




 \n {\bf Lemme 4 : }{\it Le rapport  $ \frac{z(n,p,k)}{z(q,p,k)}$ est une fonction strictement croissante de $k$ pour $k\geq 0$. }
 
  \n {\it D\émonstration  :} D'apr\ès le lemme 2, l'identit\é suivante est en effet v\érifi\ée :
  $$  \frac{z(n,p,k)}{z(q,p,k)}\equiv \frac{(q-p)!}{(n-p)!}\times\prod _{j=1}^{n-q}(q+j+k).$$
   \rightline{QED}
 

\n Soit   ${\bf k^1_p}$ le plus grand entier $h$ tel que $d\geq  \frac{z(n,p,h)}{z(q,p,h)}$.


 
 
 
 \n{\bf Th\éor\ème 6 :}
 
 {\it Supposons  le $d$-tissu de codimension $q$ ci-dessus $p$-\emph{fortement ordinaire}, $d>\frac{z(n,p,0)}{z(q,p,0)}$, et $q\geq 1$. 
 
 $(i)$ Pour $k\leq k^1_p$, l'ensemble $\widetilde{R}_k^p$ des $p$-relations ab\éliennes ferm\ées  formelles \à l'ordre $k$ poss\ède une structure naturelle de fibr\é vectoriel holomorphe de rang 
$$ \sum_{h=0}^k z(q,p,h) \Bigl(d-\frac{z(n,p,h)}{z(q,p,h)}\Bigr).$$

 $(ii)$ Le $p$-rang ferm\é $\widetilde {r}_p$ du tissu  est major\é par le nombre\footnote{Rappelons (lemme 2) la formule  $z(m,p,h)=b(m+h,m-p).c(p,h)$ .} 
 $$ \pi'_p(n,d, q):= \sum_{h\geq 0} z(q,p,h) \biggl(d-\frac{z(n,p,h)}{z(q,p,h)}\biggr)^+.$$}
 
 
 \n {\bf Remarque :} Quand il existe un $d$-tissu fortement $p$-ordinaire \emph{parall\élisable} (cf. section 6) de codimension $q$ dans un espace ambiant $n$-dimensionnel, celui-ci a un $p$-rang ferm\é maximal, et cette borne est donc optimale. 

\n {\bf Th\éor\ème 7 : }{\it Pour que le tissu soit  fortement $p$-ordinaire, il suffit que  $\widetilde P_k(p)$ soit  de rang maximum pour   $k\leq  k^1_p+1$ $($et m\^eme seulement pour $k\leq  k^1_p$ si $\widetilde  P_{k_p^1}(p)$ est une matrice carr\ée$)$.}

\n La d\émonstration des deux th\éor\èmes 6 et 7 est analogue \à celles des th\éor\èmes  3 et 4.

Toute $p$-relation ab\élienne ferm\ée est \évidemment $p$-ab\élienne, et le $p$-rang ferm\é d'un tissu est donc au plus \égal \à son  $p$-rang. Cependant, il se peut que $\pi'_p(n,d,q)$ soit plus grand que $\pi^0_p(n,d,q)$ (par exemple $\pi'_1(3,3,2)=8$ tandis que   $\pi^0_1(3,3,2)=6$). Dans ce cas, un  tissu fortement  $p$-ordinaire de rang maximal ne sera certainement pas $p$-ordinaire. Plus g\én\éralement :
\n {\bf Proposition 2 :} {\it  Pour que le tissu puisse \^etre  $p$-ordinaire, il est n\écessaire que soit r\éalis\ée la condition suivante : $$ \widetilde\alpha_k(p)-\widetilde\beta_k(p)\leq \alpha_k(p)-\beta_k(p) \hbox{ quel que soit }k\leq inf(k_p^0,k_p^1).$$}
 {\it D\émonstration : } Pour $k\leq inf(k_p^0,k_p^1)$, $\widetilde R^p_k$ doit \^etre inclus dans $ R_k^p$ d\ès que le tissu est $p$-ordinaire, et par cons\équent son rang  $\widetilde \rho_k(p)\ \bigl($toujours au moins \égal \à $\widetilde\alpha_k(p)-\widetilde\beta_k(p)\bigr)$ doit \^etre au plus \égal \à $\rho_k(p)$,  c'est-\à-dire \à $\alpha_k(p)-\beta_k(p)$ si le tissu est $p$-ordinaire. 
\section{Tissus $p$-calibr\és et connexions}

Dans cette section, nous limiterons l'expos\é au cas des  tissus calibr\és ordinaires, mais la th\éorie se transpose  sans difficult\é au cas des tissus fortement calibr\és et  fortement ordinaires.

\n Il se peut que $\frac{b(n,p).c(n,k)}{b(q,p).c(q,k)}$   prenne des valeurs enti\ères pour certaines valeurs de $k$. 

\n {\bf D\éfinition 4 :} 

\n {\it Un $d$-tissu  de codimension $q$ sera dit \emph{$p$-calibr\é} si $$d=\frac{b(n,p).c(n,k^0_p)}{b(q,p).c(q,k^0_p)}.$$ }
[si $d=\frac{z(n,p,k_p^1)}{z(q,p,k_p^1)}$ on  dit  que le tissu est  \emph{fortement} $p$-calibr\é].

\n Pour all\éger les notations, l'entier $p$ \étant bien fix\é, posons dans cette section :
 $$ R_{h}:= R^p_{h},\ \ k_0:=k^0_p\ \hbox{ et }
\ {\cal E}:=  R_{k_0-1}\ .$$
  Le fibr\é      
  $$ R_{k_0}:=J^1 {\cal E}\cap J^{k_0} R_{0}$$ est l'intersection des fibr\és $J^1 {\cal E}$ et $ J^{k_0}R_{0}$ dans $J^1 ( J^{k_0-1} R_{0} )$.

Si le tissu est $p$-calibr\é  et $p$-ordinaire, la projection $ R_{k_0} \to {\cal E}$ est  un isomorphisme de fibr\és vectoriels holomorphes 
 de rang $\pi_p^0(n,d,q)$.  
Notant  $u:{\cal E}\buildrel\cong\over \rightarrow \tilde R_{k_0}$ l'isomorphisme inverse, et 
  $\iota :\tilde R_{k_0}\subset J^1 {\cal E}$ l'inclusion naturelle, l'application compos\ée $v:=\iota \scirc u$ de $ {\cal E}$ dans $ J^1 {\cal E}$ est une scission holomorphe  de la suite exacte
  $$0\to T^*U\otimes {\cal E} \buildrel{\buildrel{\nabla}\over\longleftarrow}\over{\longrightarrow} J^1 {\cal E}\buildrel{\buildrel{v}\over\longleftarrow}\over{\longrightarrow}   {\cal E}\to 0  $$
   et d\éfinit par cons\équent  une connexion holomorphe  $\nabla $    sur ${\cal E} $, que nous appellerons   \emph{la connexion tautologique},   dont la d\érivation covariante associ\ée est donn\ée par la formule
  $$\nabla s= j^1 s -\bigl(\iota \scirc u \bigr)(s). $$

  \n {\bf Th\éor\ème 8 } :  {\it 
  
 \n $(i)$ Si un $d$-tissu de codimension $q$  est $p$-ordinaire  et  calibr\é, ses $p$-relations ab\éliennes  s'identifient, par l'application $\sigma\to j^{k_1-1}\sigma$, aux sections holomorphes $s$ de ${\cal E}$   dont la d\ériv\ée covariante  $\nabla s $ par rapport \à la connexion tautologique est nulle.
 
 \n $(ii)$ Le tissu   est alors  de rang maximum $\pi_p^0(n,d,q) $ ssi la courbure de la connexion tautologique est nulle. 
 }

 
   
    

\n {\it D\émonstration :}

  Puisque  $v$ se factorise   \à travers $ R_{k_0}$, il est \équivalent de dire, pour une section  $\sigma $ de $ R_{0}$, que $j^{k_0}  \sigma$ est une section de $ R_{k_0}$ ou que $\nabla (j^{k_0-1}  \sigma) $ s'annule :  les $p$-relations ab\éliennes  sont donc  les sections holomorphes  $\sigma $ de $ R_{0}$ telles que  $\nabla (j^{k_0-1}  \sigma)=0 $. 
  
Dire que cette  connexion  est  sans courbure \équivaut alors \à dire que  le tissu est   de rang maximum $\pi_p^0(n,d,q)$ (le rang de ${\cal E}$).

\rightline{QED}

\section{Exemples :}

Les calculs concernant ces exemples n'ont  pas toujours \ét\é d\étaill\és, quand ils \étaient   trop compliqu\és ou trop fastidieux pour \^etre faits sans ordinateur.
 \subsection{$p=q=1 :$} On retrouve  le r\ésultat de [CL] : $$\pi'_1(n,d,1) =
 \pi'(n,d)\ \ \Bigl(:=\sum_{h\geq 1} \bigl(d-c(n,h)\bigr)^+\Bigr).$$
 Le calcul de la courbure  dans le cas  calibr\é alors \ét\é pr\ésent\é dans [DL1], ainsi que de nombreux exemples 
 (voir aussi [DL2] dans le cas non-n\écessairement calibr\é).

\subsection{Tissus parall\élisables :} Un $d$-tissu de codimension $q$ est dit \emph{parall\élisable} s'il est possible de choisir
les  coordonn\ées locales et les syst\èmes g\én\érateurs  de fonctions ${\cal F}_i$-basiques $u_i=(u_{i,\alpha})_\alpha$ de fa\c con que toutes les fonctions $u_{i,\alpha}$ soient des fonctions lin\éaires de ces coordonn\ées locales. 

   \n {\bf Th\éor\ème 9 :} 
   {\it Pour des tissus parall\èlisables $p$-ordinaires $($resp. fortement $p$-ordinaires$)$, le $p$-rang $($resp. le $p$-rang ferm\é$)$ est toujours maximal, \égal \à $\pi^0_p(n,d,q)$ $\bigl($resp. \à  $\pi'_p(n,d,q)\bigr)$. Si le tissu n'est pas $p$-ordinaire $($resp. fortement $p$-ordinaire$)$,  le $p$-rang   $($resp. le $p$-rang ferm\é$)$ est strictement  plus grand. }
  
  \n {\it D\émonstration :}  Si les fonctions $u_{i,\alpha}$ sont toutes affines,   les matrices $Q_k(p)$ (resp. $\widetilde Q_k(p)$) sont alors  nulles. Le rang de ${\cal M}_k(p)$ $\bigl($resp. $\widetilde{\cal M}_k(p)\bigr)$ est donc  \égal \à la somme des rangs des $P_h(p)$ (resp. $\widetilde P_h(p)$) pour $h\leq k$. Par cons\équent, si le tissu est $p$-ordinaire (resp. fortement $p$-ordinaire), $R_p^{k_0}$ (resp. $\widetilde R_p^{k_1}$) est de rang maximum $\pi_p^0(n,d,q)$ (resp. $\pi_p'(n,d,q)$), et  le rang de  $R_p^k$ (resp. $\widetilde R_p^k$ ne diminue pas quand $k>k_p^0$ (resp. $k>k^1_p$). 
  
  \rightline{QED}

Quand on peut choisir\footnote{Toutefois, ce n'est pas toujours possible ; par exemple, si $q=n-1$, ce n'est possible que pour $p=n-1$, d'apr\ès la proposition 3 ci-dessous.} les fonctions affines $u_{i,\alpha}$ de fa\c con que les tissus soient $p$-ordinaires (resp. fortement $p$-ordinaires), 
ceci montre le caract\ère optimal des bornes des th\éor\èmes 3 et 6.

 \subsection{Tissus de courbes ($q=n-1  $) - G\én\éralit\és :}  Le nombre $\pi'_{n-1}(n,d,n-1)$ est alors \égal \à 
 $$\sum_{h=0}^{d-n-1} b(n-2+h,h)\ (d-n-h).$$On retrouve, pour les tissus en courbes $(n-1)$-ordinaires,  la borne donn\ée par Damiano ([D]). \hb Plus g\én\éralement, $$\pi'_p(n,d,n-1)=\sum_{h\geq 0} b(n-1+h,n-1-p).b(p-1+h,h)\ \biggl(d-\frac{n+h}{n-p}\biggr)^+ \hbox{pour les tissus fortement $p$-ordinaires,} $$
et $$\pi^0_p(n,d,n-1)= b(n-1,p).\sum_{h\geq 0} b(n-2+h,n-2)\ \biggl(d-\frac{n(n-1+h)}{(n-1)(n-p)}\biggr)^+\hbox{pour les tissus  $p$-ordinaires, }(1\leq p<n-1)$$

Notons $X_i$ un champ de vecteurs engendrant ${\cal F}_i$.
 \n {\bf Proposition 3 :} {\it S'il existe deux indices distincts $i,j$  tels que le crochet $[X_i,X_j]$ soit  une combinaison lin\éaire de $X_i$ et $X_j$, les p-rangs $r_p$ et $\widetilde r_p$ de celui-ci sont infinis, pour tout $p\leq n-2$. En particulier, le tissu n'est certainement ni p-ordinaire, ni fortement p-ordinaire pour ces valeurs de $p$.}
 \n  C'est un cas particulier de la proposition 1 (section 2), \évidemment sans objet pour $p=n-1$.

\subsection{Tissus de courbes en dimension 3 ($q=2$) :}
Simplifions les notations g\én\érales en d\ésignant les coordonn\ées locales par  $(x,y,z) $ au lieu de $(x_1,x_2,x_3)$, et les  fonctions ${\cal F}_i$-basiques d\éfinissant ${\cal F}_i$ par $u_i(x,y,z)$ et $v_i(x,y,z)$, au lieu de $u_{i,1}$ et $u_{i,2}$. 
\subsubsection{Simplification des expressions de $\pi'_2(3,d,2),\ \pi^0_1(3,d,2)$ et $\pi'_1(3,d,2)$ :}
Pour $n=3$, chaque terme sous le signe $\sum$ dans les expressions pr\éc\édentes est  polynomial  de degr\é 2 en $h$. En utilisant les formules classiques donnant  les sommes d'entiers et les sommes de carr\és,  on obtient    :
$$\pi'_{2}(3,d,2)=\frac{1}{6}(d-1)(d-2)(d-3) \hbox{ \ \ \ (expression d\éj\à donn\ée dans [D])},$$ 
 tandis que 
$$\pi^0_1(3,d,2)=\frac{1}{4}\delta(\delta+1)(\delta+\rho),\hbox{ o\ù l'on a pos\é  } 4d-4=3\delta+\rho , \ ( \rho=0,1,2),$$

et $$\pi'_1(3,d,2)=\frac{1}{3}(d^2-1)(2d-3).$$

\subsubsection{Exemples de 4-tissus de courbes en dimension 3 :}

D\éfinissons le 4-tissu $\cal W$, d\épendant de  deux  fonctions a priori arbitraires $\varphi(x,y,z)$ et $\psi(x,y,z) $ par  
$$u_1\equiv x,\ u_2\equiv y,\ u_3\equiv z,\ u_4\equiv \varphi(x,y,z),\ v_1\equiv y,\ v_2\equiv z,\ v_3\equiv x,\ v_4\equiv \psi(x,y,z),$$

D'apr\ès la proposition 3, ces tissus 
 ne sont certainement ni 1-ordinaires ni fortement 1-ordinaires et sont de   rang  $r_1(\cal W )$ et $\widetilde r_1(\cal W )$   infinis. Les conditions de la proposition 2 ne sont d'ailleurs pas r\éalis\ées.

Etudions leurs  2-relations ab\éliennes. 
On cherche alors des 2-formes qui peuvent s'\écrire :$$\omega_i=h_i(u_i,v_i)\ du_i\wedge dv_i.$$ 
 Notons respectivement $J_{i,xy}$, $J_{i,yz}$ et $J_{i,zx}$ les d\éterminants des matrice jacobiennes  $\frac{D(u_i,v_i)}{D(x,y)}$, $\frac{D(u_i,v_i)}{D(y,z)}$ et $\frac{D(u_i,v_i)}{D(z,x)}$.
La condition de trace nulle s'\écrit : $$<P_0(2),h>\equiv 0,\ \ \hbox{ avec }\ 
P_0(2)=\begin{pmatrix}J_{xy}\\J_{yz}\\J_{zx}\end{pmatrix},$$
$J_{xy},\ J_{yz}, \ J_{zx}$ d\ésignant les $d$-vecteurs ligne $(J_{1,xy}\cdots J_{d,xy})$, $(J_{1,yz}\cdots J_{d,yz})$ et $(J_{1,zx}\cdots J_{d,zx})$, et $h$ le $d$ vecteur colonne $(h_1,\cdots,h_d)$.

Notant respectivement $\Delta u_x$, $\Delta u_y$, $\Delta u_z$, $\Delta v_x$ $\Delta v_y$ et $\Delta v_z$ les matrices diagonales $d\times d$ construites sur les $d$-vecteurs $\bigl((u_1)'_x,\cdots,(u_d)'_x\bigr)$, $\bigl((u_1)'_y,\cdots,(u_d)'_y\bigr)$, $\bigl((u_1)'_z,\cdots,(u_d)'_z\bigr)$, $\bigl((v_1)'_x,\cdots,(v_d)'_x\bigr)$, $\bigl((v_1)'_y,\cdots,(v_d)'_y\bigr)$, et $\bigl((v_1)'_z,\cdots,(v_d)'_z\bigr)$,  on construit d'abord la matrice $9\times 2d$ par blocs 
qu'on obtient en composant avec $P_0$ chacune des 3 d\ériv\ées partielles  du $d$-vecteur $\bigl(h_1(u_1,v_1),\cdots,h_1(u_1,v_1)\bigr)$. Mais la somme des lignes 3,4 et 8 est nulle car la 2-forme $\sum_i h(u_i,v_i)\ du_i\wedge dv_i$ doit \^etre ferm\ée.
La matrice   $\widetilde P_1(2)$ (de taille $8\times 2d$) s'obtient donc en   supprimant de la matrice pr\écédente l'une de ces trois lignes, disons la huiti\ème pour fixer les id\ées :
$$\widetilde P_1(2)=\begin{pmatrix}
J_{xy}.\Delta u_x&J_{xy}.\Delta v_x\\
J_{xy}.\Delta u_y&J_{xy}.\Delta v_y\\
J_{xy}.\Delta u_z&J_{xy}.\Delta v_z\\
J_{yz}.\Delta u_x&J_{yz}.\Delta v_x\\
J_{yz}.\Delta u_y&J_{yz}.\Delta v_y\\
J_{yz}.\Delta u_z&J_{yz}.\Delta v_z\\
J_{zx}.\Delta u_x&J_{zx}.\Delta v_x\\
J_{zx}.\Delta u_z&J_{zx}.\Delta v_z\\
\end{pmatrix},
\hbox { et } \widetilde Q_1(2)=\begin{pmatrix}
(J_{xy})'_x\\
(J_{xy})'_y\\
(J_{xy})'_z\\
(J_{yz})'_x\\
(J_{yz})'_y\\
(J_{yz})'_z\\
(J_{zx})'_x\\
(J_{zx})'_z\\
\end{pmatrix}$$

 Posant \hb  $\varphi'_x=a\ ,\varphi'_y=b\ ,\varphi'_z=c $, $\psi'_x=p ,\psi'_y=q\ ,\psi'_z=r $, \ \ et \ \ $br-qc=A$, $cp-ar=B$, $aq-pb=C$, \hb on obtient :
$$ P_0(2)=\begin{pmatrix}1&0&0&C\\0&1&0&A\\0&0&1&B\end{pmatrix},
\widetilde P_1(2)=\begin{pmatrix}
a&0&0&aC&p&0&0&pC\\
0&a&0&aA&0&p&0&pA\\
0&0&a&aB&0&0&p&pB\\
b&0&0&bC&q&0&0&qC\\
0&b&0&bA&0&q&0&qA\\
0&0&b&bB&0&0&q&qB\\
c&0&0&cC&r&0&0&rC\\
0&0&c&cB&0&0&r&rB\\
\end{pmatrix},\hbox{ et }
\widetilde Q_1(2)=\begin{pmatrix}
0&0&0&C'_x\\
0&0&0&C'_y\\
0&0&0&C'_z\\
0&0&0&A'_x\\
0&0&0&A'_y\\
0&0&0&A'_z\\
0&0&0&B'_x\\
0&0&0&B'_z
\end{pmatrix}.$$
Le rang de $P_0(2)$ est toujours maximal, et  le $4$-vecteur  $$s:=(h_1=-C,h_2=-A,h_3=-B,h_4=1)$$ est une base du module des sections du fibr\é ${\cal E}=Ker P_0(2)$, de rang 1.  Puisque le d\éterminant  de $\widetilde P_1(2))$ est \égal \à $-ABC,$ le tissu est  $2$-ordinaire ssi $ABC$ n'est pas nul.

\n {\bf Cas o\ù  ${\cal F}_4$ est un feuilletage en droites parall\èles :} Les expressions  $a,b,c,p,q, r$ et $A,B,C$ sont des constantes. Les nombres $A$, $B$ et $C$ ne peuvent  pas \^etre tous nuls, puisque  $d\varphi\wedge d\psi\neq 0$.

\n - 1) {\bf si $ABC$ n'est pas nul}, le tissu est 2-ordinaire : \hb
Dans ce cas, 
$\widetilde r_2=1\ \bigl(=\pi'_2(3,d,2)\bigr)$, et $(h_1=-C,h_2=-A,h_3=-B,h_4=1)$ d\éfinit une base de $Ab^2$,  soit  
$$\omega_1=-C\ dx\wedge dy,\omega_2=-A\ dy\wedge dz, \omega_3=-B\ dz\wedge dx,\omega_4=A\ dy\wedge dz+ B\ dz\wedge dx+ C\ dx\wedge dy.$$
Ces deux relations ab\éliennes sont des cobords : la 2-relation pr\éc\édente est en effet la diff\érentielle de la 1-relation ab\élienne
$$\eta_1=-C  x\ dy,\ \eta_2=-A  y\ dz,\ \eta_3=-B  z\ dx,\ \eta_4=A y\ dz+ B z\ dx+ C  x\ dy.$$
- 2) {\bf si  l'un  des trois nombres $A,B$ ou $C$ est nul},  (disons $C$ pour fixer les id\ées), et pas les deux autres :  $h_4(ax+by+cz,px+qy+rz)$ ne doit d\épendre que de la seule variable $z$ ; autrement dit la fonction $h_4(u,v)$ doit v\érifier l'une des deux \équations \équivalentes $a h'_u+ph'_v\equiv 0$ ou $bh'_u+qh'_v\equiv 0$, soit 
$h'_u(u,v)\equiv -p \ \xi(u,v)$ et $h'_v\equiv a\ \xi(u,v)$, et le rang $\widetilde r_2$ est infini (chaque fonction $\xi(u,v)$ telle que $a\ \xi'_u+p\ \xi'_v\equiv 0$ d\éfinissant une 2-relation ab\élienne).

\n {\bf Cas g\én\éral (exemple de courbure) : } On va utiliser le fait que le tissu est fortement 2-calibr\é, et calculer -lorsqu'il est fortement 2-ordinaire- la courbure de la connexion tautologique correspondante sur le fibr\é ${\cal E}=Ker\ P_0(2)$, de rang 1 : le tissu sera de rang  $\widetilde r_2$ \égal \à 1 ou 0, selon que cette courbure est nulle ou non. 

\n Puisque la quatri\ème composante de $s$ est \égale \à 1, les seules composantes  du 8-vecteur colonne $\bigl(\widetilde P_1(2)\bigr)^{-1}.\widetilde Q_1(2).s$ qui nous int\éressent pour calculer la forme de   connexion  (relative \à la trivialisation d\éfinie par $\{s\}$) sont   la quatri\ème et la huiti\ème,  
soit :
$$H:=\frac{1}{ABC}\Bigl(pAC'_z-rA'_xC\Bigr)\hbox{ et } K:=\frac{1}{ABC}\Bigl(cA'_xC-aAC'_z\Bigr).$$
On en d\éduit la forme de connexion 
$$\eta:=H\ d\varphi+K\ d\psi \hbox{ et la forme de  courbure } \Omega:=dH\wedge d\varphi+dK\wedge d\psi.$$
Prenons par exemple la famille de tissus  ${\cal W}_\lambda$,   d\épendant d'un    param\ètre scalaire $\lambda$, obtenue avec   $$\varphi(x,y,z)\equiv x+y\ ,\hbox{ et }  \psi(x,y,z)\equiv x+z+\frac{1}{2}(x^2+2\lambda xz+z^2).$$ Ce tissu est fortement 2-ordinaire au voisinage de l'origine (en dehors des deux plans d'\équations respectives $1+\lambda x+z=0$ et $1+ x+\lambda  z=0$), et la courbure de la connexion est \égale \à 
$$\frac{\lambda(\lambda-1)(z-x)}{p^2r^2}\bigl((\lambda+1)(x+z)+2\bigr)\ dx\wedge dz .$$
Ainsi, $\widetilde r_2({\cal W}_\lambda)$ est \égal \à 1  si $\lambda =0$ ou 1, et \à 0 sinon : 

\n - {\bf si $\lambda=0$}, 
$$\omega_1=-(1+u_1)\ du_1\wedge dv_1,\ \omega_2=(1+v_2)\ du_2\wedge dv_2,\ \omega_3=-(1+u_3)\ du_3\wedge dv_3,\ \omega_4=du_4\wedge dv_4$$ d\éfinit en effet une 2-relation ab\élienne fournissant une base de $Ab^2$ ; cette 2-relation  ab\élienne  est  un cobord 
: c'est   la diff\érentielle de la 1-relation ab\élienne 
$$\eta_1=-(u_1+v_1)(1+u_1)\ du_1,\ \eta_2=-u_2(1+v_2)\ dv_2,\ \eta_3=-v_3(1+u_3)\ du_3,\ \eta_4= u_4\ dv_4.$$ 
 - {\bf si $\lambda=1$}, c'est
$$\omega_1=du_1\wedge dv_1,\ \omega_2=- du_2\wedge dv_2,\ \omega_3= du_3\wedge dv_3,\ \omega_4=(2v_4+1)^{-1/2}\ du_4\wedge dv_4$$  qui d\éfinit  une base de $Ab^2$. 
 Cette 2-relation  ab\élienne  est encore un cobord : c'est la diff\érentielle de la 1-relation ab\élienne 
$$\eta_1= (1+u_1) (du_1+dv_1),\ \eta_2=v_2\ du_2,\ \eta_3= u_3\ dv_3,\ \eta_4= - (2v_4+1)^{1/2}\ du_4.$$  



\subsection{Exemples de 4-tissus de codimension 2 en dimension 4 :}

Dans [G], V.V.Goldberg a donn\é trois exemples de 4-tissus de codimension 2 en dimension 4,  dont le 2-rang \était maximal \égal \à 1 ($=\pi_2(4,4,2)$).
Nous allons voir ci-dessous comment les distinguer par leur 1-rang et les invariants que nous avons d\éfinis. 


Notant $(x,y,z,t)$ les coordonn\ées dans {\bb C}$^4$, et $(u_i,v_i)$ au lieu de $(u_{i,1},u_{i,2})$ un syst\ème g\én\érateur de fonctions basiques des feuilletages ${\cal F}_i$ ($i=1,2,3,4$), les exemples de [G] sont les les tissus ${\cal W}_1$, ${\cal W}_2$, et ${\cal W}_3$ suivants  :

\n- pour ${\cal W}_1$ : $(u_1= x$, $v_1= y)$, $(u_2= z$, $v_2= t)$, $\bigl(u_3= x+z$, $v_3=(y+t)(z-x)\bigr)$, et\hb
\indent$\Bigl(u_4=\frac{(y+t)^2(z-x)^2}{yt}$\ ,\  $v_4=x+z+(y+t)(z-x) .\varphi\bigl((yt)^{1/2}\bigr)$,
 o\ù l'on a pos\é $\varphi(s)= \frac{{\rm Arctg} s}{s}\Bigr)$, 
\n - pour ${\cal W}_2$  : $(u_1= x$, $v_1= y)$, $(u_2= z$, $v_2= t)$, $(u_3= x+z$, $v_3= yz-xt)$, et\hb \indent $\Bigl(u_4=\frac{yz-xt}{y+t}$\ , \ $v_4=-(x+z)-\frac{yz-xt}{y+t} . ln \frac{t}{y}\Bigr)$ ;
\n- pour ${\cal W}_3$  : $(u_1= x$, $v_1= y)$, $(u_2= z$, $v_2=t)$, $(u_3= x+z+\frac{1}{2}x^2t$, $v_3= y+t-\frac{1}{2}xt^2)$, et \hb \indent$\bigl(u_4=-x+z+\frac{x^2t}{2}$, \ $v_4=y-t-\frac{xt^2}{2}\bigr)$.

\n Une base de l'espace des germes de 2-relations ab\éliennes en un point g\én\érique est donn\ée
 
\n-  par $\omega_1= \frac{1}{v_1}\ du_1\wedge dv_1\ ,\ \omega_2=\frac{1}{v_2}\ du_2\wedge dv_2\ ,\ \omega_3=-\frac{1}{v_3}\ du_3\wedge dv_3\ ,\ \omega_4= -\frac{1}{2u_4}\ du_4\wedge dv_4 $ pour ${\cal W}_1$,

\n -  par $\omega_1= -\frac{1}{v_1}\ du_1\wedge dv_1\ ,\ \omega_2=-\frac{1}{v_2}\ du_2\wedge dv_2\ ,\ \omega_3=\frac{1}{v_3}\ du_3\wedge dv_3\ ,\ \omega_4= -\frac{1}{u_4}\ du_4\wedge dv_4 $ pour ${\cal W}_2$,

\n- et par $\omega_1= 2\ du_1\wedge dv_1\ ,\ \omega_2=2\ du_2\wedge dv_2\ ,\ \omega_3=-\ du_3\wedge dv_3\ ,\ \omega_4= du_4\wedge dv_4 $ pour ${\cal W}_3$.

\n Les deux premiers exemples,  ${\cal W}_1$ et ${\cal W}_2$, admettent d'autre part une 1-relation ab\élienne ferm\ée \évidente,  donn\ée par 
$$\omega_1=-du_1,\ \omega_2=-du_2,\ \omega_3=du_3,\ \omega_4=0.$$
On peut  les distinguer  par   le rang de $P_1(1)$ qui est \égal \à 15 pour ${\cal W}_1$ et \à 13 pour ${\cal W}_2$. Ils ne sont cependant  1-ordinaires ni l'un ni l'autre, puisque le rang maximum de $P_1(1)$ est 16. D'autre part, le rang maximum de  $\widetilde P_1(1)$, \égal \à 10,  est atteint dans le cas de ${\cal W}_1$ qui est fortement 1-ordinaire, et pas dans celui de ${\cal W}_2$ pour lequel  il est \égal \à 9. 

Quant \à ${\cal W}_3$, il n'est pas non plus 1-ordinaire ni fortement 1-ordinaire, puisque les rangs des matrices $P_1(1)$ et $\widetilde P_1(1)$ sont repectivement 11 et 9.

\section{Un germe de relation ab\élienne ferm\ée qui n'est pas un cobord}

  Notons plus g\én\éralement $(x,y,z,t) $ les coordonn\ées sur  un voisinage  $U$ de l'origine dans  {\bb C}$^4$, et d\éfinissons les feuilletages ${\cal F}_i$  d'un 4-tissu ${W}_\varphi$ de codimension 2 par un syst\ème g\én\érateur  $(u_i,v_i)_i$ de fonctions :

$u_1=x,\hskip 4cm v_1=y$,  

$u_2=z,\hskip 4cmv_2=t,$ 

$ u_3=x+z+x.\varphi(x,t),\hskip 1.7cm v_3=y+t-t.\varphi(x,t)$,

$ u_4=-x+z+x.\varphi(x,t),\hskip 1.4cm v_4=y-t-t.\varphi(x,t),$

\n o\ù $\varphi$ d\ésigne une fonction holomorphe.

\n  Il existe  d'autres  0-relations ab\éliennes que celles qui sont ferm\ées,  comme le prouve la formule suivante, facile  \à v\érifier :
   
   \n {\bf Lemme 5 :} 
   $$2(u_1v_1+u_2v_2)-u_3v_3+u_4v_4=0.$$
   Ceci prouve en particulier que $r_1(W_\varphi))\geq 1$, puisque $Ab^1(W_\varphi)$ contient au moins des cobords non nuls.
   
   Notant $\rho_0=(2u_1v_1,\ 2u_2v_2,\ -u_3v_3,\ u_4v_4)$ cette relation ab\élienne, nous allons montrer le 
   
    \n {\bf Lemme 6 :}{\it
     
    \n  Pour $\varphi(x,t) =\frac{xt}{2}$, ${Ab}^1({W}_\varphi)$ est engendr\é par $d(\rho_0)$.}
    
 \n {\it D\émonstration :} 
 
 \n Soit $(\omega_i)_i$ une 1-relation ab\élienne. avec
 
\centerline{ $\omega_1=f_1(x,y) \ dx+g_1(x,y) \ dy,$\hskip 1cm $\ \omega_2=f_2(z,t) \ dz+g_2(z,t) \ dt,\ $}
$\vspace{.11cm}$
\centerline{$\omega_3=f_3(u_3,v_3) \ du_3+g_3(u_3,v_3) \ dv_3,\ \ 
\omega_4=f_4(u_4,v_4) \ du_4+g_4(u_4,v_4) \ dv_4$.}
 La condition de trace nulle s'\écrit :
$$\begin{matrix}
(f_3-f_4)+(\varphi+x\varphi'_x)(f_3+f_4)-t\varphi'_x(g_3+g_4) &=&-f_1\ ,\\\\
g_3+g_4&=&-g_1 \ ,\\\\
f_3+f_4&=&-f_2\ ,\\\\
x\varphi'_t(f_3+f_4)+(g_3-g_4)-(\varphi+x\varphi'_t)(g_3+g_4) &=&-g_2\ ,
\end{matrix} ,$$
 soit :
$$\begin{matrix} 
-2f_3&=& f_1(x,y)+(1-\alpha).f_2(z,t)+\gamma .g_1(x,y)\ ,\\  &&\\
2f_4&=& f_1(x,y)-(\alpha+1).f_2(z,t)+ \gamma. g_1(x,y)\ ,\\&&\\
-2g_3&=&  (\beta+1).g_1(x,y)+g_2(z,t)-\delta .f_2(z,t)\ ,\\&&\\
2g_4&= &(\beta-1).g_1(x,y)+g_2(z,t)-\delta .f_2(z,t) \ ,
\end{matrix}
$$
o\ù l'on a  pos\é   pour simplifier : $\alpha=\varphi+x\varphi'_x$, $\beta =\varphi+t\varphi'_t$, $\gamma=t\varphi'_x$ et $\delta=x\varphi'_t$.

Pour exprimer que $f_3$ et $g_3$ (resp. $f_4$ et $g_4$) ne sont fonctions que de $u_3$ et $v_3$ (resp. $u_4$ et $v_4$), nous allons faire le changement de coordonn\ées   $$X=u_3,\ Y=v_3,\ Z=u_4,\ T=v_4.$$
On obtient les champs de vecteurs
$$\begin{matrix} 
2\frac{\partial}{\partial X}&=&\frac{\partial}{\partial x}&+&\gamma\frac{\partial  }{\partial y}&+&(1-\alpha)\frac{\partial  }{\partial z}\ ,&&\\&&&&&&&&\\
2\frac{\partial}{\partial Y}&=&&&(1+\beta)\frac{\partial  }{\partial y}&-&\delta\frac{\partial  }{\partial z}&+&\frac{\partial  }{\partial t}\ ,\\&&&&&&&&\\
-2\frac{\partial}{\partial Z}&=&\frac{\partial}{\partial x}&+&\gamma\frac{\partial  }{\partial y}&-&(1+\alpha)\frac{\partial  }{\partial z}\ ,&\\&&&&&&&&\\
-2\frac{\partial}{\partial T}&=&&&(\beta-1)\frac{\partial  }{\partial y}&-&\delta\frac{\partial  }{\partial z}&+&\frac{\partial  }{\partial t}\ ,
\end{matrix}
$$

\n Notaént $\epsilon$ un nombre pouvant prendre la  valeur $+1$ ou $-1$,  posons :

 $F=f_1+(\epsilon-\alpha).f_2+\gamma.g_1$, \hskip 1cm $G=g_2+(\beta+\epsilon)g_1-\delta f_2$, 
   
   $D=\frac{\partial}{\partial x}+\gamma\frac{\partial  }{\partial y}-(\epsilon+\alpha)\frac{\partial  }{\partial z}\ $ \hskip 1cm
et $\widetilde{D}=(\beta-\epsilon)\frac{\partial  }{\partial y}-\delta\frac{\partial  }{\partial z}+\frac{\partial  }{\partial t}.\ $

\n Les huit  \équations 

$\frac{\partial f_3}{\partial Z}=0$, $\frac{\partial f_4}{\partial X}=0$, $\frac{\partial f_3}{\partial T}=0$, $\frac{\partial f_4}{\partial Y}=0$, $\frac{\partial g_3}{\partial Z}=0$, $\frac{\partial g_4}{\partial X}=0$, $\frac{\partial g_3}{\partial T}=0$,   $\frac{\partial g_4}{\partial Y}=0$ 

\n \équivalent \à $D(F)=0$, $\widetilde{D}(F)=0$, $D(G)=0$, $\widetilde{D}(G)=0$ avec $\epsilon=+1$ et $\epsilon=-1$.

\n En particulier, si $\varphi(x,t)=\frac{xt}{2}   $,
$$D(G)=(tg_1-xf_2+xt.A)+\epsilon A
\hbox{ avec $A=(g_1)'_x-(g_2)'_z+\frac{t^2}{2}(g_1)'_y+\frac{x^2}{2}(f_2)'_z$.}$$
Puisque cette expression doit \^etre nulle aussi bien  pour $\epsilon=+1$ que $-1$, on en d\éduit les deux identit\és
$A=0$ et $tg_1-xf_2=0$. 
Mais la seconde de ces identit\és, que l'on peut encore \écrire $\frac{g_1(x,y)}{x}=\frac{f_2(z,t)}{t}$,  ne peut \^etre r\éalis\ée  que s'il existe un scalaire $k$ tel que 
$g_1(x,y)=kx$ et $f_2(z,t)=kt$. 
En reportant ces expressions dans l'\équation $A=0$,  on obtient : $(g_2)'_z=k$, et en reportant ces r\ésultats  dans les \équations 
 $\widetilde D(F)=0$ et $\widetilde D(G)=0$, on obtient : $f_1(x,y)=ky$ et $g_2(z,t)=kz$.
Puisque $f_3,f_4,g_3$ et $g_4$ sont d\étermin\és en fonction de $f_1,f_2,g_1$ et $g_2$, on en d\éduit que le 1-rang du tissu est au plus 1. On savait d\éja qu'il \était au moins 1, puisque $d(\rho_0)$ est une relation ab\élienne, d'o\ù la conclusion du lemme.   On peut aussi calculer directement 
$$f_3=-\frac{1}{2}kv_3,\ g_3=-\frac{1}{2}ku_3, \ f_4=kv_4 \hbox{ et }g_4=ku_4.$$

\rightline{QED}
D'autre part, si  $\varphi(x,t) =\frac{xt}{2},$ Goldberg a d\émontr\é dans $[G]$ que le 2-rang  $r_2({W}_\varphi)$ \était \égal \à 1, $Ab^2({W}_\varphi)$ \étant engendr\é par la relation ab\élienne $$(2du_1\wedge dv_1,\ 2du_2\wedge dv_2 ,\ -du_3\wedge dv_3,\ du_4\wedge dv_4).$$ On en d\éduit :

\n {\bf Proposition 4 :}

{\it

$r_1({W}_\varphi)\geq 1$

\n Pour $\varphi(x,t) =\frac{xt}{2},$ 

$r_1({W}_\varphi)=1$, et 
$ Ab^1({W}_\varphi))$ ne contient que des cobords,

  $H^2_{Ab}({W}_\varphi)\cong\hbox{\bb C}$.
}


\begin{thebibliography}{dango 9999}


\bibitem[A]{A} M.A. Akivis, {\it Differential geometry of webs}, J. Soviet Math. 29, Springer, 1985,1631-1647.

\bibitem[B]{B} W. Blaschke et G. Bol, {\it Geometrie der Gewebe}, Die Grundlehren der Mathematik 49, Springer, 1938.

\bibitem[Bo]{Bo} G. Bol, {\it \"Uber ein bemerkenswertes F\"unfgewebe in der Ebene}, Abh. Math. Hamburg Univ., 11, 1936, 387-393.

\bibitem[C]{C} S.S. Chern, {\it Abz\"ahlungen f\"ur Gewebe,} Abh. Math. Hamb. Univ. 11, 1935, 163-170.

\bibitem[CG]{CG}  S.S. Chern,  P.A. Griffiths, {\it An inequality for the rank of a web, and webs of maximum rank}, Ann. Scuola Norm. Sup. Pisa 5, 1978, 539-557.

\bibitem[CL]{CL} V. Cavalier, D. Lehmann, {\it Ordinary holomorphic webs of codimension one. } arXiv 0703596v2 [mathsDS], 2007, et Ann. Sc. Norm. Super. Pisa, cl. Sci (5), vol XI (2012), 197-214.

\bibitem[D]{D} D.B. Damiano, {\it Abelian equations and characteristic classes}, Thesis, Brown University, (1980) ;  American J. Math. 105-6, 
1983,  1325-1345.


 
\bibitem[DL1]{DL1} J. P. Dufour, D. Lehmann, {\it Calcul explicite de la courbure des tissus calibr\'es ordinaires,; } arXiv 1408.3909v1 [mathsDG], 18/08/2014.

\bibitem[DL2]{DL2} J. P. Dufour, D. Lehmann, {\it Rank of ordinary webs in codimension one : an effective method ; } arXiv 1703.03725v1 [math.DG], 10/03/2017.



\bibitem[G]{G}  V.V. Goldberg,  {\it Theory of multi-codimensional webs}, Kluwer, Dordrecht, 1988.

\bibitem[GHL]{GHL} L. Gruson, Y. Hantout, D. Lehmann, {\it  Courbes alg\ébriques ordinaires et tissus associ\és},  Comptes Rendus de l'Acad\émie des Sciences, Paris, s\ér. I, 350 (2012), 513-518.

\bibitem[Gr]{Gr}  P.A. Griffiths, {\it On Abel differential equations}, Algebraic Geometry, The Johns Hopkins centennial Lectures, Ed. J.-I. Igusa (1977), 26-51.

\bibitem[H1]{H1} A. H\'enaut, {\it Planar web geometry through abelian relations  and connections}
   Annals of Math. 159 (2004)  425-445.
   
 


 \bibitem[H2] {H2} A. H\énaut, {\it Formes diff\érentielles ab\éliennes, bornes de Castelnuovo et g\éom\étrie des tissus}, Commentarii Math.  Helvetici, 79 (1), 2004, 25-57.







  
 

\bibitem[Pa]{Pa} A. Pantazi, Sur la d\étermination du rang d'un tissu plan, C.R. Acad. Sc. Roumanie 2, 1938, 108-111.

\bibitem[P]{P} L. Pirio, {\it Equations Fonctionnelles Ab\'eliennes et G\'eom\'etrie des tissus},
  Th\`ese de doctorat de l'Universit\'e Paris VI, 2004.
  
  \bigskip
  
  Daniel Lehmann,  \hb  lehm.dan@gmail.com
  
  ancien professeur \à l'Universit\é de Montpellier 2,
  
  4 rue Becagrun, 30980 Saint Dionisy, France.









\end{thebibliography}
\end{document}